\documentclass[]{article}
\input cyracc.def

\usepackage{amssymb,amsmath,amsthm,bm,citesort}

\newcommand{\R}{\mathbf{R}}

\renewcommand{\P}{\mathrm{P}}
\newcommand {\E}{\mathrm{E}}

\renewcommand{\d}{\text{\rm d}}
\newcommand{\e}{\text{\rm e}}

\newcommand{\lip}{\mathrm{Lip}}

\newtheorem{stat}{Statement}[section]
\newtheorem{proposition}[stat]{Proposition}

\newtheorem{theorem}[stat]{Theorem}
\newtheorem{lemma}[stat]{Lemma}
\theoremstyle{definition}
\newtheorem{definition}[stat]{Definition}\newtheorem{remark}[stat]{Remark}

\numberwithin{equation}{section}

\begin{document}

\title{\bf On the global maximum of
	the solution to a stochastic heat equation with
	compact-support initial data\thanks{%
	Research supported in part by NSF grant DMS-0704024.}}
	
\author{Mohammud Foondun \and Davar Khoshnevisan}

\date{January 19, 2009}
\maketitle
\begin{abstract}
Consider a stochastic heat equation
$\partial_t u = \kappa  \partial^2_{xx}u+\sigma(u)\dot{w}$
for a space-time white noise $\dot{w}$ and a
constant $\kappa>0$. Under some suitable conditions on the initial function $u_0$ and $\sigma$, we show that the quantity
\begin{equation*}
	\limsup_{t\to\infty}t^{-1}\ln\E\left(\sup_{x\in\R}
	|u_t(x)|^2\right)
\end{equation*}
is bounded away from zero and infinity by explicit multiples of $1/\kappa$.
Our proof works by demonstrating quantitatively that the peaks of the
stochastic process $x\mapsto u_t(x)$
are highly concentrated for infinitely-many
large values of $t$.
In the special case of the parabolic Anderson model---where $\sigma(u)=
\lambda u$ for some $\lambda>0$---this ``peaking'' is a way to
make precise the notion of physical intermittency.

		
		
\noindent{\it Running Title:} The stochastic heat equation.
\end{abstract}

\section{Introduction}

We consider the stochastic heat equation,
\begin{equation}\label{heat:cont}
	\frac{\partial u_t(x)}{\partial t} =
	\kappa\frac{\partial^2 u_t(x)}{\partial x^2} +
	\sigma(u_t(x))\dot{w}(t\,,x)
	\quad\text{for $t>0$ and $x\in\R$},
\end{equation}
where $\kappa>0$ is fixed,
$\sigma:\R\to\R$ is Lipschitz continuous with $\sigma(0)=0$,
$\dot{w}$ denotes space-time white noise, and the initial data
$u_0:\R\to\R$ is nonrandom. There are several areas to which
\eqref{heat:cont} has deep and natural connections;
perhaps chief among them are the
stochastic Burgers' equation \cite{GN} and the celebrated
KPZ equation of statistical mechanics \cite{Kardar,KPZ};
see also \cite[Chapter 9]{KrugSpohn}.

It is well known that \eqref{heat:cont}
has an almost-surely unique, adapted and continuous
solution $\{u_t(x)\}_{t\ge 0,x\in\R}$ \cite[Theorem 6.4, p.\ 26]{SPDEBook}.
In addition, the condition that $\sigma(0)=0$ implies
that if $u_0\in L^2(\R)$, then
$u_t\in L^2(\R)$ a.s.\ for all $t\ge 0$; see Dalang and Mueller
\cite{DM}. Note that our conditions on $\sigma$ ensure that
\begin{equation}
	|\sigma(u)|\le\lip_\sigma|u|\qquad\text{for all $u\in\R$},
\end{equation}
where
\begin{equation}
	\lip_\sigma:=\sup_{-\infty<x<x'<\infty}
	\left| \frac{\sigma(x)-\sigma(x')}{x-x'}\right|.
\end{equation}
Our goal is to establish the following
general growth estimate.

\begin{theorem}\label{th:max:LE}
	Suppose there exists ${\rm L}_\sigma\in(0\,,\infty)$
	such that $|\sigma(u)|\ge{\rm L}_\sigma|u|$ for all $u\in \R$.
	Suppose also that $u_0\not\equiv 0$ is H\"older-continuous
	of order $\ge 1/2$,
	nonnegative, and supported
	in $[-K\,,K]$ for some finite $K>0$.
	Then, \eqref{heat:cont} has an almost-surely unique,
	continuous and adapted solution $\{u_t(x)\}_{t\ge 0,x\in\R}$ such that
	$u_t\in L^2(\R)$ a.s.\ for all $t\ge 0$, and
	\begin{equation}\label{eq:max:LE}
		\frac{{\rm L}^4_\sigma}{8\kappa}\le
		\limsup_{t\to\infty} t^{-1}
		\ln \E\left( \sup_{x\in\R}|u_t(x)|^2\right)
		\le\frac{\lip^4_\sigma}{8\kappa}.
	\end{equation}
\end{theorem}

Because of Mueller's comparison principle \cite{Mueller},
the nonnegativity of $u_0$ implies
that $\sup_{t,x}\E(|u_t(x)|)=\sup_{t,x}
\E(u_t(x))$, and this quantity has to be finite because
$u_0$ is bounded; confer with \eqref{mild}.
Consequently,
\begin{equation}\label{eq:ll}
		\sup_{x\in\R} \left\| u_t(x)\right\|_{L^1(\P)} \ll
		\left\|\sup_{x\in\R} u_t(x)\right\|_{L^2(\P)}
		\quad\text{as $t\to\infty$}.
\end{equation}

When $\lip_\sigma
= {\rm L}_\sigma$, \eqref{heat:cont} becomes
the well-studied parabolic Anderson model
\cite{BC,CM94}. And \eqref{eq:ll} makes
precise the physical notion that the solution to \eqref{heat:cont}
concentrates near ``very high peaks''
\cite{BC,CM94,Kardar,KPZ}.

In order to explain the idea behind our proof, we introduce the following.

\begin{definition}\label{def:eff-comp}
	We say that a continuous random field $f:=\{f(t\,,x)\}_{t\ge 0,x\in\R}$ 
	has \emph{effectively-compact support}
	if there exists a nonrandom measurable function $p:\R_+\to\R_+$ of at-most polynomial
	growth such that:
	\begin{enumerate}
		\item[(a)]  $\limsup_{t\rightarrow \infty}
			t^{-1}\ln \int_{|x|\le  p(t)}\E(|f(t\,,x)|^2)\,\d x >0$; and
		\item[(b)]$\limsup_{t\to\infty} t^{-1}
			\ln\int_{|x|> p(t)}\E(|f(t\,,x)|^2)\,\d x <0$.
	\end{enumerate}
	We might refer to the function $p$ as the \emph{radius of effective support} of $f$.
	\qed
\end{definition}
One of the ideas here is to use Mueller's comparison principle
\cite{Mueller} to compare $\sup_{x\in\R}|u_t(x)|$ with 
the $L^2(\R)$-norm of
$x\mapsto u_t(x)$, which is easier to analyze.
We carry these steps out in Lemma \ref{lem:nonunif}. 
We also appeal to the fact that
the compact-support property of $u_0$ implies that 
$u_t(x)$ has an effectively-compact support
[Proposition \ref{prop:u-eff-comp}].  This can be interpreted as 
a kind of optimal regularity theorem. However, these matters need to be handled delicately,
as ``effectively compact'' cannot be replaced
by ``compact''; see Mueller \cite{Mueller}.

Our method for establishing an effectively-compact
support property is motivated strongly by ideas of Mueller and
Perkins \cite{MuellerPerkins}. In the cases that $u_t(x)$
denotes the density of some particles at $x$ at time $t$,
our effectively-compact support property implies that most of
the particles accumulate on a very small set. This method
might appeal to the
reader who is interested in mathematical
descriptions of physical intermittency.

Throughout this paper we use the mild formulation of
the solution, in accordance with Walsh \cite{Walsh}.
That is, $u$ is the a.s.-unique adapted solution to
\begin{equation}\label{mild}
	u_t(x)=(p_t*u_0)(x) + \int_0^t\int_{-\infty}^\infty
	p_{t-s}(y-x)\sigma(u_s(y))\,w(\d s\,\d y),
\end{equation}
where $p_\tau(z):=(4\kappa\tau\pi)^{-1/2}
\exp(-z^2/(4\kappa\tau))$ denotes
the heat kernel corresponding to the operator
$\kappa\partial^2/\partial x^2$,
and the stochastic integral is understood in the sense of
Walsh \cite{Walsh}. Some times we write $\|X\|_p$ in place of
$\{\E(|X|^p)\}^{1/p}$.

\section{A preliminary result}

As mentioned in the introduction, the strategy behind our proof of Theorem \ref{th:max:LE}
is to relate the global maximum of the solution
to a ``closed-form quantity'' that resembles $\sup_x|u_t(x)|$
for large values of $t$. That closed-form quantity turns out to be
the $L^2(\R)$-norm of $x\mapsto u_t(x)$. Our next result
analyses the growth of the mentioned closed-form quantity.
We related it to $\sup_x|u_t(x)|$ in the next section.
The methods of this section follow closely the classical ideas of
Choquet and Deny \cite{CD}
that were developed in a determinstic setting.

\begin{theorem}\label{th:L2}
	Suppose $\sigma:\R\to\R$ is Lipschitz continuous,
	$\sigma(0)=0$, and there exists ${\rm L}_\sigma\in(0\,,\infty)$
	such that ${\rm L}_\sigma|u|\le |\sigma(u)|$ for all $u\in \R$.
	If $u_0\in L^2(\R)$ and $u_0\not\equiv 0$,
	then \eqref{heat:cont} has an almost-surely unique,
	continuous and adapted solution $\{u_t(x)\}_{t\ge 0,x\in\R}$ such that
	$u_t\in L^2(\R)$ a.s.\ for all $t\ge 0$, and
	\begin{equation}\label{heat:cont:LE}
		\frac{{\rm L}_\sigma^4}{8\kappa}\le
		\limsup_{t\to\infty} t^{-1}
		\ln \E\left(\|u_t\|_{L^2(\R)}^2\right)
		\le \frac{\lip_\sigma^4}{8\kappa}.
	\end{equation}
\end{theorem}

\begin{proof}
	It suffices to establish \eqref{heat:cont:LE}.  Note that
	\begin{equation}\begin{split}
	&\E\left( \left| u_t(x)\right|^2\right)\\
		&\qquad=\left| (p_t* u_0)(x)\right|^2 +\int_0^t
		\d s\int_{-\infty}^\infty\d y\
		\E\left(\left| \sigma(u_s(y))\right|^2\right)
		\cdot\left| p_{t-s}(y-x)\right|^2.
	\end{split}\end{equation}

	Because $|\sigma(u)|\ge {\rm L}_\sigma |u|$,
	\begin{equation}\label{eq:pre:planch}\begin{split}
		&\E\left(\| u_t\|_{L^2(\R)}^2\right)\\
		&\qquad=\|p_t*u_0\|_{L^2(\R)}^2
			+\int_0^t\d s\int_{-\infty}^\infty\d y\
			\E\left(\left| \sigma(u_s(y))\right|^2\right)
			\cdot\| p_{t-s}\|_{L^2(\R)}^2\\
		&\qquad\ge \|p_t*u_0\|_{L^2(\R)}^2+{\rm L}_\sigma^2\cdot
			\int_0^t\E\left(\| u_s\|_{L^2(\R)}^2\right)
			\cdot\| p_{t-s}\|_{L^2(\R)}^2\, \d s.
	\end{split}\end{equation}
	We can multiply the preceding by $\exp(-\lambda t)$ throughout
	and integrate $[\d t]$ to find that if
	\begin{equation}\label{def:L}
		 U(\lambda):=\int_0^\infty \e^{-\lambda t}
		\E\left(\| u_t\|_{L^2(\R)}^2\right)\, \d t,
	\end{equation}
	then
	\begin{equation}
		 U(\lambda) \ge \int_0^\infty \e^{-\lambda t}
		\|p_t*u_0\|_{L^2(\R)}^2\,\d t +{\rm L}_\sigma^2\cdot
		 U(\lambda)\cdot\int_0^\infty \e^{-\lambda t}
		\|p_t\|_{L^2(\R)}^2\,\d t.
	\end{equation}
	According to Plancherel's theorem, the following holds for all
	finite Borel measures $\mu$ on $\R$:
	\begin{equation}
		\|p_t*\mu\|_{L^2(\R)}^2 = \frac{1}{2\pi}
		\int_{-\infty}^\infty \left| \hat\mu(\xi)
		\right|^2\e^{-2\kappa t\xi^2}\,\d\xi.
	\end{equation}
	Therefore, Tonelli's theorem ensures that
	\begin{equation}
		\int_{0}^\infty \e^{-\lambda t}	\|p_t*\mu\|_{L^2(\R)}^2
		\,\d t = \frac{1}{2\pi}\int_{0}^\infty
		\frac{|\hat\mu(\xi)|^2}{\lambda+2\kappa\xi^2}\,\d\xi.
	\end{equation}
	We apply this identity twice in \eqref{eq:pre:planch}:
	Once with $\mu:=\delta_0$; and once with
	$\d\mu/\d x :=u_0$. This leads us to the following.
	\begin{equation}\label{eq:L2:ineq1}\begin{split}
		 U(\lambda) &\ge \frac{1}{2\pi}\int_{0}^\infty
			\frac{|\hat{u}_0(\xi)|^2}{\lambda+
			2\kappa\xi^2}\,\d\xi
			+ {\rm L}_\sigma^2 \cdot U(\lambda)
			\cdot\frac{1}{2\pi}\int_{0}^\infty
			\frac{1}{\lambda+2\kappa\xi^2}\,\d\xi\\
		&=\frac{1}{2\pi}\int_{0}^\infty
			\frac{|\hat{u}_0(\xi)|^2}{\lambda+2\kappa\xi^2}\,\d\xi
			+ {\rm L}_\sigma^2 \cdot U(\lambda)
			\cdot\frac{1}{2\sqrt{2\kappa\lambda}}.
	\end{split}\end{equation}
	Since $u_0\not\equiv 0$, the first [Fourier]
	integral is strictly positive.
	Consequently, the above recursive relation shows that $ U(\lambda) =\infty$
	if $\lambda\le  {\rm L}_\sigma^4/(8\kappa)$.
	This and a real-variable argument
	together imply the first inequality in \eqref{heat:cont:LE}; see \cite{FK} for more details.
	
	For the other bound we use a Picard-iteration argument in order
	to obtain an {\it a priori} estimate. Let $u^{(0)}_t(x):=u_0(x)$ and iteratively define
	\begin{equation}
		u^{(n+1)}_t(x) := (p_t*u_0)(x) + \int_0^t\int_{-\infty}^\infty 
		p_{t-s}(y-x)\sigma( u^{(n)}_s(y))\, w (\d s\,\d y).
	\end{equation}
	Since $\|p_t*u_0\|_{L^2(\R)}\le
	\|u_0\|_{L^2(\R)}$ and $|\sigma(u)|\le\lip_\sigma|u|$,
	\begin{equation}\label{eq:L2:ineq2}\begin{split}
		&\E\left(\left\| u_t^{(n+1)}\right\|_{L^2(\R)}^2\right)\\
		&\hskip.5in\le \|u_0\|_{L^2(\R)}^2+ \lip_\sigma^2\cdot
			\int_0^t\E\left(\left\| u_s^{(n)}\right\|_{L^2(\R)}^2\right)
			\cdot\| p_{t-s}\|_{L^2(\R)}^2\, \d s.
	\end{split}\end{equation}
	Therefore, if we set
	\begin{equation}
		 M^{(k)} (\lambda):= \sup_{t\ge 0}\left[\e^{-\lambda t}
		\E\left(\left\|u_t^{(k)}\right\|_{L^2(\R)}^2\right)\right],
	\end{equation}
	then it follows that
	\begin{equation}\begin{split}
		 M^{(n+1)} (\lambda)
			&\le \|u_0 \|_{L^2(\R)}^2 +\lip_\sigma^2\cdot
			 M^{(n)} (\lambda)\cdot \int_0^\infty \e^{-\lambda(t-s)}
			\left\| p_{t-s}\right\|_{L^2(\R)}^2\,\d s\\
		&= \|u_0 \|_{L^2(\R)}^2
			+\frac{\lip_\sigma^2}{2\sqrt{2\kappa\lambda}}\,
			 M^{(n)} (\lambda).
	\end{split}\end{equation}
	Thus, in particular, $\sup_{n\ge 0} M^{(n)} (\lambda)<\infty$
	if $\lambda> \lip_\sigma^4/(8\kappa)$. We can argue similarly to
	show also that if $\lambda>\lip_\sigma^4/(8\kappa)$, then
	\begin{equation}
		\sum_n\sup_{t\ge 0}\left[ \e^{-\lambda t}\E
		\left(\left\| u^{(n+1)}_t-u^{(n)}_t\right\|_{L^2(\R)}^2 \right)
		\right]^{1/2}<\infty.
	\end{equation}
	In particular, uniqueness shows that if $\lambda>\lip_\sigma^4/(8\kappa)$, then
	\begin{equation}
		\lim_{n\to\infty} \sup_{t\ge 0}\left[ \e^{-\lambda t}\E
		\left(\left\| u^{(n)}_t-u_t\right\|_{L^2(\R)}^2 \right)
		\right]=0.
	\end{equation}
	Consequently,  if $\lambda>\lip_\sigma^4/(8\kappa)$, then
	\begin{equation}\begin{split}
		\sup_{t\ge 0}\left[\e^{-\lambda t}\E\left(\|u_t\|_{L^2(\R)}^2\right)\right]
			&=\lim_{n\to\infty} M^{(n)}(\lambda)\\
		&\le\sup_{k\ge 0} M^{(k)}(\lambda)<\infty.
	\end{split}\end{equation}
	The second inequality of
	\eqref{heat:cont:LE} follows readily from this bound.
\end{proof}

\section{Proof of Theorem \ref{th:max:LE}}

Our proof of Theorem \ref{th:max:LE} hinges on a number of
steps, which we develop separately.
First we recall the following.

\begin{proposition}[Theorem 2.1 and Example 2.9 of \cite{FK}]\label{pr:FK}
	If $u_0$ is bounded and measurable, then $u_t(x)\in L^p(\P)$
	for all $p\in[1\,,\infty)$. Moreover, $\overline\gamma(p)<\infty$ for all
	$p\in[1\,,\infty)$ and $\overline{\gamma}(2)\le\lip_\sigma^4/
	(8\kappa)$, where
	\begin{equation}
		\overline{\gamma}(p):= \limsup_{t\to\infty}
		t^{-1}  \sup_{x\in\R} \ln
		\E\left(|u_t(x)|^p\right).
	\end{equation}
\end{proposition}

Next, we record a simple though crucial property of
the function $\overline\gamma$.

\begin{remark}\label{rem:cvx}
	Suppose $X$ is a nonnegative random variable
	with finite moments of all order. By H\"older's inequality,
	$p\mapsto \ln \E(X^p)$ is convex on $[1\,,\infty)$.
	It follows that $\overline\gamma$ is convex---in particular
	continuous---on $[1\,,\infty)$.\qed
\end{remark}

Now we begin our analysis, in earnest, by deriving an
upper bound on the $L^k(\P)$-norm of the solution
$u_t(x)$ that includes simultaneously
a sharp decay rate in $x$ and a sharp explosion rate in $t$.

\begin{lemma}\label{lem:nonunif}
	Suppose that $u_0\not\equiv 0$, and $u_0$ is supported
	in $[-K,K]$ for some finite constant $K>0$.
	Then, for all real numbers $k\in[1\,,\infty)$ and $p\in(1\,,\infty)$,
	\begin{equation}\label{eq:nonunif}
		\limsup_{t\to\infty}  t^{-1} \sup_{x\in\R}
		\left( \frac{x^2}{4t^2}+\frac{k+1-(1/p)}{k}\ln \E\left(|u_t(x)|^k\right)
		\right)\le \frac{\overline\gamma(kp)}{p}.
	\end{equation}
\end{lemma}

\begin{proof}
	According to Mueller's comparison principle
	(\cite{Mueller}; more specifically, see \cite[Theorem 5.1, p.\ 130]{SPDEBook}),
	the solution to \eqref{heat:cont}
	has the following nonnegativity property: Because $u_0\ge 0$
	then outside a single null set,
	$u_t\ge 0$ for all $t\ge 0$. And therefore,
	\begin{equation}
		\left\| u_t(x)\right\|_1 =
		(p_t*u_0)(x)=\frac{1}{\sqrt{4\kappa \pi t}}
		\int_{-K}^K \e^{-(x-y)^2/(4\kappa t)} u_0(y)\,\d y.
	\end{equation}
	Because $(x-y)^2\ge (x^2/2)-K^2$,
	\begin{equation}\label{eq:L1:bd}
		\|u_t(x)\|_1 \le \text{const}\cdot
		\e^{-x^2/(4t)}\qquad\text{for all $x\in\R$ and $t\ge 1$}.
	\end{equation}
	The constant appearing in the above display depends on $K$.
	Next we note that for every $\theta\in(0\,,\infty)$,
	\begin{equation}\label{eq:Lk:bd}\begin{split}
		\E\left(|u_t(x)|^k\right) &\le \theta^k + \E\left(
			|u_t(x)|^k\,; u_t(x)\ge\theta\right)\\
		&\le \theta^k + \left( \E\left(|u_t(x)|^{kp}\right)
			\right)^{1/p}\cdot\left( \P\{u_t(x)>\theta \}\right)^{1-(1/p)}.
	\end{split}\end{equation}
	Proposition \ref{pr:FK} implies that
	\begin{equation}\label{eq:Lkp:bd}
		\sup_{x\in\R}\left( \E\left(|u_t(x)|^{kp}\right) \right)^{1/p}
		\le \exp\left( t\ \frac{\overline\gamma(kp)+o(1)}{p}
		\right),
	\end{equation}
	where $o(1)\to 0$ as $t\to\infty$. Also, we can apply \eqref{eq:L1:bd}
	together with the Chebyshev inequality to find that
	\begin{equation}\label{eq:cheb:bd}
		\left( \P\left\{ u_t(x)>\theta\right\}\right)^{1-(1/p)} \le
		\text{const}\cdot \theta^{-1+(1/p)}
		\exp\left( -\frac{x^2}{4t}\cdot\left[
		1-\frac1p\right]\right).
	\end{equation}
	Taking into consideration \eqref{eq:Lkp:bd} and \eqref{eq:cheb:bd},
	inequality \eqref{eq:Lk:bd} reduces to
	\begin{equation}\begin{split}
		\E\left(|u_t(x)|^k\right) &\le \inf_{\theta>0}
			\left( \theta^k + \alpha \theta^{-1+(1/p)}\right),\\
	\end{split}\end{equation}
	where
	\begin{equation}
		\alpha := \exp\left( -\frac{x^2}{4t}\cdot\left[
		1-\frac1p\right]+t\ \frac{\overline\gamma(kp)+o(1)}{p}\right).
	\end{equation}
Some calculus shows that the function 
$g(\theta):=(\theta^k+\alpha\theta^{-1+1/p})1_{(\theta>0)}$
attains its minimum at $\theta:= ( (p-1)/kp )^{p/(kp+p-1)}$. This yields
\begin{equation*}
	\E\left(|u_t(x)|^k\right) \le  \alpha^{kp/(kp+p-1)}\left(\frac{p-1}{kp}
	\right)^{kp/(kp+p-1)}\cdot\left( \frac{1-p-kp}{1-p}
	\right).
\end{equation*}
We now divide both sides of the above display by
$\alpha^{kp/(kp+p-1)}$ and take the appropriate limit to obtain the result.
\end{proof}

Our next lemma is a basic estimate of continuity in the
variable $x$. It is not entirely standard as it holds uniformly
for all times $t\ge 0$. We emphasize that the constant $p$
is assumed to be an integer.
We will deal with this shortcoming subsequently.

\begin{lemma}\label{lem:modulus}
	Suppose that the initial function $u_0$ is H\"older continuous of order $\geq 1/2$.
	Then, for all integers $p\ge 1$
	and $\beta>\overline\gamma(2p)$
	there exists a constant $A_{p,\beta}\in(0\,,\infty)$ such that
	the following holds: Simultaneously for all $t\ge 0$,
	\begin{equation}\label{eq:modulus}
		\sup_{j\in\mathbf{Z}}\sup_{j\le x<x'\le j+1}
		\left\| \frac{ u_t(x)-u_t(x') }{|x-x'|^{1/2}}\right\|_{2p} \le
		A_{p,\beta}\, \e^{\beta t/(2p)}.
	\end{equation}
\end{lemma}

\begin{proof}
	Burkholder's inequality \cite{Burkholder} and Minkowski's inequality together
	imply that
	\begin{equation}\label{eq:pre-est}\begin{split}
		&\left\| u_t(x)-u_t(x')\right\|_{2p}\\
		&\le \left| (p_t*u_0)(x)- (p_t*u_0)(x')\right|\\
		& \quad+ z_{2p}
			\left\| \int_0^t \d s\int_{-\infty}^\infty \d y\
			\left|\sigma(u_s(y))\right|^2\cdot
			\left| p_{t-s}(y-x)-p_{t-s}(y-x')\right|^2
			\right\|_p^{1/2}\\
		&\le \left| (p_t*u_0)(x)- (p_t*u_0)(x')\right|\\
		& \quad+ z_{2p}'
			\left\| \int_0^t \d s\int_{-\infty}^\infty \d y\
			\left| u_s(y)\right|^2\cdot
			\left| p_{t-s}(y-x)-p_{t-s}(y-x')\right|^2
			\right\|_p^{1/2},
	\end{split}\end{equation}
	where $z_p$ is a positive and finite constant
	that depend only on $p$, and $z_p':=z_p \lip_\sigma$.
	
	On one hand, we have the following consequence
	of Young's inequality:
	\begin{equation}\label{eq:T1}\begin{split}
		\sup_{t\ge 0}\sup_{|x-x'|\le\delta}
			\left| (p_t*u_0)(x)- (p_t*u_0)(x')\right|&\le
			\sup_{|a-b|\le\delta}\left| u_0(a)-u_0(b)\right|\\
		&\le \text{const}\cdot\delta^{1/2}.
	\end{split}\end{equation}
	
	On the other hand, the generalized H\"older inequality suggests
	that if $p\ge 1$ is an integer, then
	for all $s_1,\ldots,s_p\ge 0$ and $y_1,\ldots,y_p\in\R$,
	\begin{equation}\label{eq:GHolder}
		\E\left(\prod_{j=1}^p \left| u_{s_j}(y_j)\right|^p\right)
		\le \prod_{j=1}^p\left\| u_{s_j}(y_j)\right\|_{2p}^2.
	\end{equation}
	Therefore,
	\begin{equation}\begin{split}
		&\left\| \int_0^t \d s\int_{-\infty}^\infty \d y\
			\left| u_s(y)\right|^2\cdot
			\left| p_{t-s}(y-x)-p_{t-s}(y-x')\right|^2
			\right\|_p\\
		&\hskip.7in\le \int_0^t\d s\int_{-\infty}^\infty\d y\
			\|u_s(y)\|_{2p}^2\cdot \left|
			p_{t-s}(y-x)-p_{t-s}(y-x')\right|^2.
	\end{split}\end{equation}
	[Write the $p$-th power of the left-hand side as
	the expectation of a product and apply \eqref{eq:GHolder}.]
	
	A proof by contradiction shows that Proposition
	\ref{pr:FK} gives the following [see \cite{FK} for more details]:
	\begin{equation}
		c_\beta := \sup_{s\ge 0}\sup_{y\in\R} \left[\e^{-\beta s}
		\E\left( |u_s(y) |^{2p}\right)\right] <\infty
		\qquad\text{for all $\beta>\overline\gamma(2p)$}.
	\end{equation}
	Consequently,
	\begin{equation}\begin{split}
		&\left\| \int_0^t \d s\int_{-\infty}^\infty \d y\
			\left| u_s(y)\right|^2\cdot
			\left| p_{t-s}(y-x)-p_{t-s}(y-x')\right|^2
			\right\|_p\\
		&\hskip.5in\le c_\beta^{1/p}\cdot
			\int_0^t\d s\int_{-\infty}^\infty\d y\
			\e^{\beta s/p}\cdot \left|
			p_{t-s}(y-x)-p_{t-s}(y-x')\right|^2\\
		&\hskip.5in\le c_\beta^{1/p}\e^{\beta t/p}\cdot
			\int_0^\infty \d s\ \e^{-\beta s/p}
			\int_{-\infty}^\infty\d y\ \left|
			p_s(y-x)-p_s(y-x')\right|^2.
	\end{split}\end{equation}
	Since $\hat{p}_s(\xi)=\exp(-\kappa s\xi^2)$, Plancherel's theorem
	tells us that the right-hand side of the preceding inequality is equal to
	\begin{equation}\begin{split}
		\frac{c_\beta^{1/p}\e^{\beta t/p}}{\pi}\cdot
			\int_0^\infty \d s\ \e^{-\beta s/p}
			&\int_{-\infty}^\infty\d \xi\ \e^{-2\kappa
			s\xi^2}\left[
			1-\cos(\xi(x-x'))\right]\\
		&\ = \frac{2c_\beta^{1/p}\e^{\beta t/p}}{\pi}\cdot
			\int_0^\infty\frac{[1-\cos(\xi(x-x'))]}{(\beta/p)+
			2\kappa\xi^2} \,\d\xi.
	\end{split}\end{equation}
	Because $1-\cos\theta\le \min(1\,,\theta^2)$, a direct estimation of
	the integral leads to the following bound:
	\begin{equation}\begin{split}
		&\left\| \int_0^t \d s\int_{-\infty}^\infty \d y\
			\left| u_s(y)\right|^2\cdot
			\left| p_{t-s}(y-x)-p_{t-s}(y-x')\right|^2
			\right\|_p\\
		&\hskip2.6in\le \text{const}\cdot \e^{\beta t/p} \cdot
			|x-x'|,
	\end{split}\end{equation}
	where the implied constant depends only on $p$,
	$\kappa$, and $\beta$.
	This, \eqref{eq:T1}, and \eqref{eq:pre-est} together imply the lemma.
\end{proof}

The preceding lemma holds for all {\it integers} $p\ge 1$. In the
following, we improve it [at a slight cost] to the case that
$p\in(1\,,2)$ is a real number.

\begin{lemma}\label{lem:unif-modulus}
	Suppose the conditions of Lemma \ref{lem:modulus} are met.
	Then for all $p\in(1\,,2)$ and $\delta\in(0\,,1)$
	there exists a constant $B_{p,\delta}\in(0\,,\infty)$ such that
	the following holds: Simultaneously for all $t\ge 0$
	and $x,x'\in\R$ with $|x-x'|\le 1$,
	\begin{equation}\label{eq:unif-modulus}
		\E\left( |u_t(x)-u_t(x')|^{2p}\right) \le
		B_{p,\delta} \cdot |x-x'|^p\cdot \e^{(1+\delta)\lambda_pt},
	\end{equation}
	where
	\begin{equation}\label{eq:lambda}
		\lambda_p := (2-p)\overline\gamma(2)+
		(p-1)\overline\gamma(4).
	\end{equation}
\end{lemma}

\begin{proof}
	We start by writing
	\begin{equation*}
		\E\left(|u_t(x)-u_t(x')|^{2p}\right)=
		\E\left(|u_t(x)-u_t(x')|^{2(2-p)}|u_t(x)-u_t(x')|^{4(p-1)}\right).
	\end{equation*}

	We can apply H\"older's inequality to conclude that for all
		$p\in(1\,,2)$, $t\ge 0$, and $x,x'\in\R$,
		\begin{equation}\begin{split}
			&\E\left(|u_t(x)-u_t(x')|^{2p}\right) \\
			&\hskip.8in\le \left[ \E\left( |u_t(x)-u_t(x')|^2\right) \right]^{2-p}
				\left[ \E\left(|u_t(x)-u_t(x')|^4\right)\right]^{p-1}.
		\end{split}\end{equation}
	We now use Lemma \ref{lem:modulus} to obtain the following:
	\begin{equation*}
		\left[\E\left( |u_t(x)-u_t(x')|^2\right)\right]^{2-p}\le  |x-x'|^{(2-p)}A_{1,
		\beta_1}^{2(2-p)}e^{\beta_1(2-p) t}
	\end{equation*}
	and
	\begin{equation*}
		\left[ \E\left( |u_t(x)-u_t(x')|^4\right) \right]^{p-1}\le  |x-x'|^{2(p-1)}A_{2,
		\beta_2}^{4(p-1)}e^{\beta_2(p-1) t},
	\end{equation*}
	where $A_{1,\beta_1}\,,A_{2,\beta_2}\in(0,\infty)$
	and $\beta_1>\bar{\gamma}(2)$ and 
	$\beta_2>\bar{\gamma}(4)$ are fixed and finite
	constants. The proof now 
	follows by combining the above and choosing 
	$\beta_1$ and $\beta_2$ such that 
	$(1+\delta)\bar{\gamma}(2)>\beta_1>\bar{\gamma}(2)$ 
	and $(1+\delta)\bar{\gamma}(4)>\beta_2>\bar{\gamma}(4)$.
\end{proof}

The preceding lemma allows for a uniform modulus of continuity
estimate, which we record next.

\begin{lemma}\label{lem:kol}
	Suppose the conditions of Lemma \ref{lem:modulus} are met.
	Then for all $p\in(1\,,2)$ and $\epsilon,\delta\in(0\,,1)$
	there exists $C_{p,\epsilon,\delta}\in(0\,,\infty)$ such
	that simultaneously for all $t\ge 0$,
	\begin{equation}\label{eq:kol}
		\sup_{j\in\mathbf{Z}}
		\left\| \sup_{j\le x<x'\le j+1}
		\frac{|u_t(x)-u_t(x')|^2}{|x-x'|^{1-\epsilon}}\right\|_p\\
		\le C_{p,\epsilon,\delta}\cdot\e^{(1+\delta)
		\lambda_pt},
	\end{equation}
	where $\lambda_p$ was defined in \eqref{eq:lambda}.
\end{lemma}

\begin{proof}

The proof consists of an application of the Kolmogorov continuity 
theorem.  Recall that the spatial dimension is 1. 
Since $p>1$ in Lemma \ref{lem:unif-modulus}, 
we can use a suitable version of Kolmogorov continuity theorem, 
for example Theorem 4.3
of reference \cite[p.\ 10]{SPDEBook}, to obtain the result. 
The stated dependence of the constant, $C_{p,\epsilon,\delta}$ 
is consequence of the explicit form of inequality \eqref{eq:unif-modulus} 
and the proof of Theorem 4.3 in \cite{SPDEBook}.\end{proof}

Before we begin our proof of Theorem \ref{th:max:LE}, we prove that under some condition the $L^2(\P)$-norm of the solution has an effectively-compact support.
\begin{proposition}\label{prop:u-eff-comp}
	If the conditions of Theorem \ref{th:max:LE} are met, then there exists
	a finite and positive constant $m$ such that
	$u_t(x)$ has an effectively-compact support
	with radius of effective support $p(t)=mt$.
\end{proposition}

\begin{proof}
	We begin by noting that for all $m,t>0$,
	\begin{equation}
		\int_{|x|>mt} |u_t(x)|^2\,\d x
		\le \int_{|x|>mt} u_t(x)\,\d x
		+\int_{\substack{|x|>mt\\
		u_t(x)\ge 1}} |u_t(x)|^2\,\d x.
	\end{equation}
	Therefore,
	\begin{equation}\label{eq:eff-comp}\begin{split}
		&\E\left(\int_{|x|>mt} |u_t(x)|^2\,\d x \right)\\
		&\le \int_{|x|>mt} (p_t*u_0)(x)\,\d x
			+\int_{|x|>mt} \E\left(|u_t(x)|^2;\, u_t(x)\ge 1\right)\,\d x.
	\end{split}\end{equation}
	Since $u_0$ has compact support,
	\eqref{eq:L1:bd} implies that
	\begin{equation}
		\int_{|x|>mt}(p_t*u_0)(x)\,\d x = O\left(  \e^{-m^2t/2}\right)\qquad
		\text{as $t\to\infty$}.
	\end{equation}
	Next we estimate the final integral in \eqref{eq:eff-comp}.
	
	Thanks to \eqref{eq:L1:bd} and Chebyshev's inequality,
	\begin{equation}
		\P\left\{ u_t(x) \ge 1\right\} \le
		\text{const}\cdot \e^{-x^2/(4t)},
	\end{equation}
	uniformly for all $x\in \R$ and $t\geq 1$. Also, from
	Proposition \ref{pr:FK}, there exists a constant $b\in(0\,,\infty)$ such that
	\begin{equation}\label{eq:LE4}
		\sup_{x\in\R}\E\left( |u_t(x)|^4\right) \le b \e^{bt/4}
		\qquad\text{for all $t\ge 1$}.
	\end{equation}
	Using the preceding two inequalities, the right-hand side 
	of inequality \eqref{eq:eff-comp} reduces to
	\begin{equation}\begin{split}
		&\E\left(\int_{|x|>mt} |u_t(x)|^2\,\d x \right)\\
		&\le  O\left( \e^{-m^2t/2}\right)+\text{const}\cdot
			\int_{|x|>mt}
			\sqrt{\E\left( |u_t(x)|^4\right)}\,\e^{-x^2/(8t)}\,\d x\\
		&\le  O\left( \e^{-m^2t/2}\right)+ \text{const}\cdot
			b^{1/2}\e^{bt/8}\cdot
			\int_{|x|>mt} \e^{-x^2/(8t)}\,\d x.
	\end{split}
	\end{equation}
	We now choose and fix $m>\sqrt{b}$ to obtain from the preceding that
	\begin{equation}
		\limsup_{t\to\infty}
		t^{-1}\ln \E\left(\int_{|x|>mt} |u_t(x)|^2\,\d x \right) <0.
	\end{equation}
	This implies part (b) of Definition \ref{def:eff-comp} with
	$p(t)=mt$. We now prove the remaining part of 
	Definition \ref{def:eff-comp}. From Theorem \ref{th:L2} and the 
	preceding, we obtain 
	for infinitely-many values of $t\rightarrow \infty$:
	\begin{equation}\label{eq:lower-comp}\begin{split}
		\exp\left( \left[\frac{{\rm L}_\sigma^4}{8\kappa}+o(1)\right]t
			\right) &\le \E\left(\int_{-\infty}^\infty |u_t(x)|^2\,\d x\right)\\
		&= \E\left(\int_{-mt}^{mt} |u_t(x)|^2\,\d x\right) +
			o(1).\\
	\end{split}\end{equation}
	This finishes the proof.
\end{proof}

We will need the following elementary real-variable lemma from the theory of
slowly-varying functions. It is without doubt well known; we include a
derivation for the sake of completeness only.

\begin{lemma}\label{lem:rv}
	For every $q,\eta\in(0\,,\infty)$,
	\begin{equation}
		\int_{\e}^\infty \exp\left( -\frac{q(\ln x)^{\eta+1}}{t}\right)\,\d x
		=O\left( t^{1/\eta}\exp\left\{(t/q)^{1/\eta}\right\}
		\right)\qquad\text{as $t\to\infty$.}
	\end{equation}
\end{lemma}

\begin{proof}
	The proof uses some standard tricks. First we write the integral as
	\begin{equation}
		\int_{\e}^\infty \e^{-q(\ln x)^{\eta+1}/t}\,\d x
		= \int_1^\infty \e^{-q z^{\eta+1}/t}\, \e^z\, \d z.
	\end{equation}
	Next we change variables $[w:=z/\theta]$, for an arbitrary
	$\theta>0$, and find that
	\begin{equation}
		\int_{\e}^\infty \e^{-q(\ln x)^{\eta+1}/t}\,\d x
		=\theta\int_{1/\theta}^\infty \exp\left( -\frac{q\theta^{\eta+1}}{t} w^{\eta+1}
		+\theta w\right)\,\d w.
	\end{equation}
	Upon choosing $\theta:=(t/q)^{1/\eta}$, we obtain
	\begin{equation*}
		-\frac{q\theta^{\eta+1}}{t} w^{\eta+1}+\theta w=
		\left(\frac{t}{q} \right)^{1/\eta}(w-w^{\eta+1}),
	\end{equation*}
	and this yields
	\begin{equation}
		\int_{\e}^\infty \e^{-q(\ln x)^{\eta+1}/t}\,\d x
		= (t/q)^{1/\eta}\int_{(q/t)^{1/\eta}}^\infty
		\e^{(t/q)^{1/\eta} \cdot (w-w^{\eta+1})}\,\d w.
	\end{equation}
	Therefore, for $t$ sufficiently large, we split the
	integral on the right-hand side of the previous display as follows:
	\begin{equation}
		\int_{\e}^\infty \e^{-q(\ln x)^{\eta+1}/t}\,\d x
		= (t/q)^{1/\eta}\left( I_1+I_2\right),
	\end{equation}
	where
	\begin{equation}\begin{split}	
		I_1 &:=\int_{(q/t)^{1/\eta}}^1
			\exp\left( (t/q)^{1/\eta} \cdot (w-w^{\eta+1})
			\right)\,\d w,\\
		I_2 &:=\int_1^\infty
			\exp\left( -(t/q)^{1/\eta} \cdot w(w^\eta-1)
			\right)\,\d w.
	\end{split}\end{equation}
	Clearly,
	\begin{equation}
		I_2 \le 1+\int_2^\infty\exp\left( -(2^\eta-1)(t/q)^{1/\eta} \cdot w
		\right)\,\d w =O(1).
	\end{equation}
	The lemma follows because the
	integrand of $I_1$ is at most $\exp((t/q)^{1/\eta})$.
\end{proof}
We are now ready to establish Theorem \ref{th:max:LE}.

\begin{proof}[Proof of Theorem \ref{th:max:LE}]
	The proof of the first inequality in \eqref{eq:max:LE} is a 
	continuation of the proof Proposition \ref{prop:u-eff-comp}. 
	Indeed, from \eqref{eq:lower-comp}, we obtain
	\begin{equation}\begin{split}
		\exp\left( \left[\frac{{\rm L}_\sigma^4}{8\kappa}+o(1)\right]t
			\right) &\le \E\left(\int_{-\infty}^\infty |u_t(x)|^2\,\d x\right)\\
		&\le 2mt\cdot \E\left(\sup_{x\in\R}|u_t(x)|^2\right) +o(1).
	\end{split}\end{equation}
	We obtain first inequality in \eqref{eq:max:LE} after taking the appropriate limit.
	
	Next we prove the second inequality in \eqref{eq:max:LE} by first
	observing that for every $j\ge 1$, all increasing sequence of real numbers
	$\{a_j\}_{j=1}^\infty$ with $\sup_{j\ge 1}(a_{j+1}-a_j)\le 1$,
	$p\in(1\,,2)$, $\epsilon\in(0\,,1)$, and  $t\ge 0$,
	\begin{equation}\begin{split}
		\sup_{a_j\le x\le a_{j+1}}|u_t(x)|^{2p}
			&=\sup_{a_j\le x\le a_{j+1}}|u_t(a_j)+u_t(x)-u_t(a_j)|^{2p}\\
		&\le 2^{2p-1}
			\left( |u_t(a_j)|^{2p}
			+ \sup_{a_j\le x\le a_{j+1}} \left| u_t(x)-u_t(a_j)\right|^{2p}\right)\\
		&\le2^{2p-1}\left( |u_t(a_j)|^{2p}
			+ \left(a_{j+1}-a_j\right)^{p(1-\epsilon)}\Omega_j^p
			\right),
	\end{split}\end{equation}
	where
	\begin{equation}
		\Omega_j := \sup_{a_j\le x<x'\le a_{j+1}}
		\frac{|u_t(x)-u_t(x')|^2}{|x-x'|^{1-\epsilon}}.
	\end{equation}
	Consequently,
	\begin{equation}\begin{split}
		&\E\left(\sup_{a_j\le x\le a_{j+1}}|u_t(x)|^{2p} \right)\\
		&\hskip1in\le2^{2p-1}\left(\E\left( |u_t(a_j)|^{2p}\right)
			+ \left(a_{j+1}-a_j\right)^{p(1-\epsilon)}\E\left(\Omega_j^p\right)
			\right).
	\end{split}\end{equation}
	We use inequality \eqref{eq:nonunif} of Lemma \ref{lem:nonunif} with $k:=2p$ and $x:=a_j$ to find that
	\begin{equation}
		\E\left(|u_t(a_j)|^{2p}\right) \le\text{const}\cdot
		\exp\left( \beta_p\cdot\left[ t\ \frac{\overline
		\gamma(2p^2)+o(1)}{p}
		- \frac{a_j^2}{4t^2} \right] \right),
	\end{equation}
	where
	\begin{equation}
		\beta_p:= \frac{p}{p+1-(1/p)},
	\end{equation}
	the implied constant does not depend on $j$ or $t$,
	and $o(1)\to 0$ as $t\to\infty$, uniformly for all $j$.	
	Also, Lemma \ref{lem:kol} implies that
	\begin{equation}
		\sup_{j\ge 1}
		\E\left(\Omega_j^p\right) \le  C_{p,\epsilon,\delta}\cdot\e^{p(1+\delta)
		\lambda_pt},
	\end{equation}
	where $\delta$ is an arbitrarily-small positive constant, which
	we will choose and fix appropriately later on.
	We can combine the preceding inequalities to deduce that
	\begin{equation}\begin{split}
		&\E\left(\sup_{a_j\le x\le a_{j+1}}|u_t(x)|^{2p} \right)\\
		&\hskip1in\le \text{const}\cdot \e^{-\beta_p a_j^2/(4t^2)}\cdot
			\e^{\beta_p t (\overline\gamma(2p^2)+o(1))/p}\\
		&\hskip2in+\text{const}\cdot (a_{j+1}-a_j)^{p(1-\epsilon)}
			\e^{p(1+\delta)\lambda_pt}.
	\end{split}\end{equation}
	Choose and fix an integer $\nu\ge 1$.
	We apply the preceding with $p(1-\epsilon)>1$; we also choose the $a_l$'s
	so that $a_1:=0$, $0\le a_{j+1}-a_j\le 1$ for all $j\ge 1$,
	and $a_j:=(\log j)^\nu$ for all $j$ sufficiently large.
	Because $a_{j+1}-a_j=O((\ln j)^\nu/j)$ as $j\to\infty$,
	\begin{equation}
		\sum_{j=1}^\infty \left( a_{j+1}-a_j\right)^{p(1-\epsilon)}<\infty.
	\end{equation}
	Also, for all $J>1+\e$, sufficiently large,
	\begin{equation}\begin{split}
		\sum_{j=J}^\infty \e^{-\beta_p a_j^2/(4t^2)} &\le\int_{J-1}^\infty
			\e^{-\beta_p (\ln x)^{2\nu}/(4t^2)}\,\d x\\
		&=O\left(t^{2/(2\nu-1)} \e^{(4t^2/\beta_p)^{1/(2\nu-1)}}\right)
			\qquad(t\to\infty),
	\end{split}\end{equation}
	where we have used Lemma \ref{lem:rv}  for the last equality. 
	We can choose $\nu:=\frac12(\delta^{-1}+1)$
	so that $1/(2\nu-1)=\delta$. We can combine these terms to deduce the following:
	\begin{equation}\begin{split}
		\E\left(\sup_{x\ge a_J}|u_t(x)|^{2p}\right)&\le \sum_{j=J}^\infty
			\E\left(\sup_{a_j\le x\le a_{j+1}}|u_t(x)|^{2p} \right)\\
		&=O\left(
			t^{2\delta}\e^{(4t^2/\beta_p)^\delta+ \beta_p t (\overline\gamma(2p^2)+o(1))/p}
			+ \e^{p(1+\delta)\lambda_pt} \right).
	\end{split}\end{equation}
	A similar---though slightly simpler---argument yields the bound
	\begin{equation}
		\E\left(\sup_{0\le x\le a_J}|u_t(x)|^{2p}\right)
		= O \left(
		t^{2\delta} \e^{(4t^2/\beta_p)^\delta+ \beta_p t (\overline\gamma(2p^2)+o(1))/p}
		+ \e^{p(1+\delta)\lambda_pt} \right).
	\end{equation}
	We now use symmetry and let $\delta \downarrow 0$,
	\begin{equation}
		\limsup_{t\to\infty}t^{-1}\ln
		\E\left(\sup_{x\in\R}|u_t(x)|^{2p}\right)
		\le \max\left\{\frac{\beta_p\overline\gamma(2p^2)}{p}
		~,~ p\lambda_p \right\}.
	\end{equation}
	Let us substitute the evaluation of
	$\beta_p$ in terms of $p$ to find that
	\begin{equation}
		\limsup_{t\to\infty}t^{-1}\ln
		\E\left(\sup_{x\in\R}|u_t(x)|^{2p}\right)
		\le \max\left\{\frac{\overline\gamma(2p^2)}{p+1-(1/p)}
		~,~ p\lambda_p \right\}.
	\end{equation}
	This and Jensen's inequality together prove that
	\begin{equation}
		\limsup_{t\to\infty}t^{-1}\ln
		\E\left(\sup_{x\in\R}|u_t(x)|^2\right)
		\le \frac1p
		\max\left\{\frac{\overline\gamma(2p^2)}{p+1-(1/p)}
		~,~ p\lambda_p \right\},
	\end{equation}
	and this valid for all $p\in(1\,,2)$. As $p\downarrow 1$,
	$\lambda_p\to\overline\gamma(2)$. Moreover,
	$\overline\gamma(2p^2)\to\overline\gamma(2)$ because
	$\overline\gamma$ is convex and hence continuous on $[1\,,\infty)$
	[Remark \ref{rem:cvx}]. It follows that
	\begin{equation}
		\limsup_{t\to\infty}t^{-1}\ln
		\E\left(\sup_{x\in\R}|u_t(x)|^2\right)
		\le \overline\gamma(2),
	\end{equation}
	and this is $\le\lip_\sigma^4/(8\kappa)$
	by Proposition \ref{pr:FK}.
	\end{proof}
\begin{small}\vskip.4cm

\noindent\textbf{Mohammud Foondun \& Davar Khoshnevisan}\\
\noindent Department of Mathematics, University of Utah,
		Salt Lake City, UT 84112-0090\\
\noindent\emph{Emails:} \texttt{mohammud@math.utah.edu} \&
	\texttt{davar@math.utah.edu}\\
\noindent\emph{URLs:} \texttt{http://www.math.utah.edu/\~{}mohammud} \&
	\texttt{http://www.math.utah.edu/\~{}davar}

\end{small}

\begin{thebibliography}{99}
%
\bibitem{BC} Bertini, Lorenzo and Nicoletta Cancrini (1995).
	The stochastic heat equation: Feynman--Kac formula and
	intermittence, {\it J. Statist.\ Physics} {\bf 78}{\it (5/6)},
	1377--1402.
%
\bibitem{Burkholder} Burkholder, D. L. (1973).
	Distribution function inequalities for martingales,
	{\it Ann.\ Probab.} {\bf 1}, 19--42.
%
\bibitem{CM94} Carmona, Ren\'e A. and S. A. Molchanov (1994).
	Parabolic Anderson Problem and Intermittency,
	\emph{Memoires of the AMS} {\bf 108},
	Amer.\ Math.\ Soc., Rhode Island.
%
\bibitem{CD} Choquet, G. and J. Deny (1960).
	Sur l'\'equation de convolution
	$\mu=\mu*\sigma$, {\it C. R. Acad.\ Sci.\ Paris}
	{\bf 250}, 799--801.
%
\bibitem{SPDEBook} Dalang, Robert, Davar Khoshnevisan, Carl Mueller,
	David Nualart, and Yimin Xiao (2009).
	\emph{A Minicourse on Stochastic Partial Differential Equations},
	Lecture Notes in Mathematics {\bf 1962}
	(Davar Khoshnevisan and Firas Rassoul-Agha, editors),
	Springer, Berlin.
%
\bibitem{DM} Dalang, Robert C., and Carl Mueller (2003).
	Some non-linear s.p.d.e.'s that are second order in time,
	{\it Electron.\ J. Probab.}, Vol.\ 8, Paper no.\ 1, 1--21
	(electronic).
%
\bibitem{FKN} Foondun, Mohammud, Davar Khoshnevisan  and Eulalia Nualart(2008).
	A local time correspondence for stochastic partial differential equations.
	(Preprint)
%
\bibitem{FK} Foondun, Mohammud and Davar Khoshnevisan (2008).
	Intermittency and nonlinear parabolic stochastic partial
	differential equations.
	(Preprint)
%
\bibitem{GN} Gy\"ongy, Istv\'an and David Nualart (1999).
	On the stochastic Burgers' equation in the real line,
	{\it Ann.\ Probab.} {\bf 27}{\it (2)}, 782--802.
%
\bibitem{Kardar} Kardar, Mehran (1987).
	Replica Bethe ansatz studies of two-dimensional interfaces
	with quenched random impurities,
	\emph{Nuclear Phys.} {\bf B290}, 582--602.
%
\bibitem{KPZ} Kardar, Mehran, Giorgio Parisi, and Yi-Cheng
	Zhang (1986). Dynamic scaling of growing interfaces,
	\emph{Phys.\ Rev.\ Lett.} \textbf{56}{\it (9)}, 889--892.
%
\bibitem{KrugSpohn} Krug, J. and H. Spohn (1991).
	Kinetic roughening of growing surfaces,
	in: \emph{Solids Far From Equilibrium: Growth,
	Morphology, and Defects} (C. Godr\`eche, editor),
	pp.\ 479--582, Cambridge University Press, Cambridge.
%
\bibitem{Mueller} Mueller, Carl (1991).
	On the support of solutions to the heat equation with noise,
	\emph{Stochastics and Stoch.\ Reports}\
	\textbf{37}{\it (4)}, 225--245.
%
\bibitem{MuellerPerkins} Mueller, Carl and Edwin A. Perkins (1992).
	The compact support property for solutions to the heat equation
	with noise, {\it Probab.\ Th.\ Rel.\ Fields} {\bf 93}{\it (3)},
	325--358.
%
\bibitem{Walsh} Walsh, John B. (1986).
	\emph{An Introduction to Stochastic Partial Differential Equations},
	in: \'Ecole d'\'et\'e de probabilit\'es de Saint-Flour XIV,
	1984, pp.\ 265--439,
	Lecture Notes in Math.\ {\bf 1180}, Springer, Berlin.
%
\end{thebibliography}
\end{document}